\documentstyle[12pt]{article}
\newcommand{\be}{\begin{equation}}
      \newcommand{\ee}{\end{equation}}
      \newcommand{\ba}{\begin{eqnarray}}
       \newcommand{\ea}{\end{eqnarray}}
\newcommand{\ban}{\begin{eqnarray*}}
       \newcommand{\ean}{\end{eqnarray*}}

\newcommand{\pt}{\partial}

 \renewcommand{\o}[2]{\frac{#1}{#2}}
\newcommand{\hf}{\o{1}{2}}

\newcommand{\qed}{\hspace*{\fill}\rule{3mm}{3mm}\quad}
 \newcommand{\Pf}{\noindent {\em Proof.} }

\newcommand{\Rk}{\noindent {\em Remark.} }
\newcommand{\Def}{\noindent {\em Definition.} }
\newcommand{\bg}{\bar{g}}

\newcommand{\sign}{\mbox{\rm sign}}

\newcommand{\Tr}{\mbox{\rm Tr}}

\newcommand{\re}{\mbox{\rm Re}}

\newcommand{\n}{\nabla}

\newcommand{\bn}{\bar{\nabla}}

\font\BBb=msbm10 at 12pt 
\newcommand{\Bbb}[1]{\mbox{\BBb #1}}

\newtheorem{theo}{Theorem}[section]

\begin{document}
\newtheorem{lem}[theo]{Lemma}
\newtheorem{prop}[theo]{Proposition}  
\newtheorem{coro}[theo]{Corollary}

\title{ Eta Invariant, Conformally Flat Structures and Hyperbolic Structures}
\author{Xianzhe Dai}

\maketitle
\begin{abstract} We define an abelian group, the conformal cobordism group of hyperbolic structures, which classifies the hyperbolic structures according to whether it bounds a (higher dimensional) conformally flat structure in a conformally invariant way. We then construct a homomorphism from this group to the circle group, using the eta invariant. The homomorphism can be highly nontrivial. It remains an interesting question of how to compute this group.
\end{abstract}

\section{Introduction}

This note is inspired by a question from my colleague Darren Long (see their very interesting paper \cite{lr}). Motivated in part by considerations in physics, Long and Reid considered the question of whether a closed orientable hyperbolic $3$-manifold can be the totally geodesic boundary of a compact hyperbolic $4$-manifold. The eta invariant is shown to be an obstruction
\cite{lr}. Similarly, they considered the question of whether every flat $3$-manifold is a cusp cross-section of a complete finite volume $1$-cusped hyperbolic $4$-manifolds and showed that again the eta invariant is an obstruction.

In this note we study the question in the conformal category and show that there are more structures here. More precisely, we first define a notion of conformal cobordism among conformally flat structures according to whether one bounds a (higher dimensional) conformally flat structure in a conformally invariant way. The suitable boundary restriction here turns out to be that of the total umbilicity, which is the conformally invariant analogue of total geodesicity (Cf. Lemma \ref{tu}). Again, the eta invariant gives us an obstruction for a conformally flat structure to be null cobordant in the conformal category in that it will have to be integers. 

Hyperbolic manifolds are of course conformally flat. Note also that a cusp cross-section of a complete finite volume $1$-cusped hyperbolic manifolds is totally umbilic. Thus our result puts  the two results of \cite{lr} mentioned above into a uniform framework.

These consideration leads to the definition of an abelian group, the conformal cobordism group of hyperbolic structures, which classifies the hyperbolic structures according to whether it bounds a conformally flat structure in a conformally invariant way. The eta invariant then gives rise to a homomorphism from this conformal cobordism  group to the circle group. The homomorphism can be highly nontrivial. For example, in dimension three it has dense image. It remains a very interesting question of how to compute this group.

\section{Conformally flat structures}

Let $M$ be a smooth manifold of dimension $n$. A conformal structure on $M$ is specified by an equivalence class of Riemannian metrics:
\[ [g]=\{ \bg \ \  | \bg = e^{2\phi} g, \ \ \phi \in C^{\infty}(M) \}. \]
A conformal structure $[g]$ is called conformally flat if, for any $x\in M$, there is a coordinate neighborhood $U$ of $x$, such that 
\[ g|_U =\lambda ( \sum dx_i^2 ) , \]
where $\lambda$ is a (positive) function on $U$. We will denote by ${\cal C\! F}(n)$ the set of all conformally flat structures $(M, [g])$ of dimension $n$. 

We would be interested in conformally invariant notions about submanifolds. Consider a submanifold $N\subset M$. Its second fundamental forms with respect to the metrics $\bg$, $g$ will be denoted by $\bar{A}$ and $A$ respectively. 

\begin{prop} For any normal vector field $v$ of $N$ the bilinear forms $g(\bar{A}, v)$ and $g(A, v)$ differ by a scalar. Hence the multiplicities of the eigenvalues of the second fundamental form of $N\subset M$ are conformal invariants. In particular the total umbilicity of a hypersurface is a conformally invariant notion.
\end{prop}

Recall that a submanifold is called totally umbilic if all eigenvalues of its second fundamental form are equal. 

\Pf Let $\bg = e^{2\phi} g $ and $\bn$, $\n$ the Levi-Civita connection of $\bg$, $g$, respectively. Then their difference is a $(2, 1)$-tensor given by
\[ S(X, Y)= \bn_X Y - \n_X Y = (X\phi)Y + (Y\phi)X - g(X, Y) \n \phi. \]
Thus for a submanifold $N\subset M$, its second fundamental forms with respect to the metrics $\bg$, $g$ are related by
\be \label{cdosff} \bar{A}(X, Y) = A(X, Y) - g(X, Y) (\n \phi)_n , \ee
where 
\[ \n \phi= (\n \phi)_t + (\n \phi)_n\]
is the decomposition into the tangential and normal components. It follows then
\[ g(\bar{A}(X, Y), v)= g(A(X, Y), v) -  g(X, Y) g((\n \phi)_n, v) \]
differs only by a scalar. \qed

This leads us to the following notion of conformal cobordism. We will restrict ourself to compact oriented manifolds for the rest of the note.
\newline

\Def We say that two conformally flat structures $(M_1, [g_1])$, $(M_2, [g_2])$ are conformally cobordant if there is a conformally flat manifold $(W, [h])$ such that \\
\noindent 1). $\pt W = M_1 \cup (-M_2)$ and $[h|_{M_i}]= [g_i ]$. \\
\noindent 2).  $M_i$'s are totally umbilic submanifolds of $W$. 
\newline

Here $-M_2$ denotes $M_2$ with the opposite orientation. We point out however that, unlike the usual cobordism, it is not clear that our conformal cobordism defines an equivalent relation. But, it does if we restrict ourself to the subclass of hyperbolic structures. This observation enables us to define the conformal cobordism group of hyperbolic structures.
\newline 

\Def In the set of all hyperbolic structures of dimension $n$, the conformal cobordism defines an equivalence relation. Let $\Omega_n^{CF}(H)$ be the set of classes of conformally cobordant hyperbolic structures of dimension $n$. It is an abelian group with the group operation given by the disjoint union. We will call $\Omega_n^{CF}(H)$ the conformal cobordism group of hyperbolic structures. Its zero element is given by the collection of hyperbolic manifolds which are totally umbilic boundaries of conformally flat manifolds.
\newline

In the definition above we used the fact that if $M$ is hyperbolic, then $([0, 1] \times M, [h])$ with $h=dx^2 + g$ the product metric is conformally flat.

We also consider the connected sum of conformally flat structures. If $M_1$, $M_2$ are conformally flat manifolds of dimension $n$ and $D_j\subset M_j$ are round disks, let $M_{0j}
=M_j- ({\rm int} D_j)$. Since one can invert in round sphere it is clear that one can glue $M_{01}$, $M_{02}$ along their boundaries so that $M_1 \# M_2$ also admits a conformally flat structure.

We end this section with the following lemma which will be used later.

\begin{lem} \label{tu}
A two-sided hypersurface is totally umbilic if and only if it is totally geodesic with respect to a conformally equivalent metric. 
\end{lem}

\Pf If a hypersurface is totally geodesic with respect to a conformally equivalent metric, then the eigenvalues of its second fundamental form with respect to this conformally equivalent metric  are all zero. Hence it is totally umbilic by the conformal invariance of the multiplicity.

Conversely, if $N\subset (M, g)$ is two-sided and totally umbilic, we construct a conformal deformation of $g$ such that $N$ becomes totally geodesic. Since $N$ is two-sided we can construct  a smooth function $r$ which coincides with the signed distance to $N$ within a collar neighborhood of $N$. We use $y$ to denote (local) coordinates of $N$. Together with $r$, $(r, y)$ forms a local coordinates of $M$ in the collar neighborhood of $N$. By the total umbilicity of $N$, the eigenvalues of its second fundamental form $A$ are all equal to $\lambda(y)$ for some smooth function $\lambda(y)$ on $N$. Now set
\[ \phi = \lambda(y) r \]
in the collar neighborhood of $N$ and extend to all of $M$ smoothly in any manner. Then $N$ will be totally geodesic with respect to $\bg=e^{2\phi} g$. Indeed, according to (\ref{cdosff}), we just need to compute the normal component of $\n \phi$:
\[ (\n \phi)_n =  \lambda(y) \o{\pt}{\pt r}. \]
Our claim follows. \qed

\section{The eta invariant}

The eta invariant was introduced by Atiyah-Patodi-Singer in their seminal work on index formula for manifolds with boundary \cite{aps}. For an odd dimensional Riemannian manifold $(M^n, g)$, the signature operator is defined to be
\be \label{sigop}
A=\tau(d + \delta) = d\tau + \tau d
\ee
where $\tau = (\sqrt{-1})^{\o{n+1}{2} + p(p-1)} *$ on $p$-forms and is essentially the Hodge star operator. The signature operator depends on the metric through the Hodge star operator. One has $\tau^2=1$. 

\Rk The signature operator $A$ could be further restricted to act on the space of even (resp. odd) forms, with the two restrictions related in terms of conjugation by $\tau$. Thus our $A$ essentially decomposes into two copies of another operator, which is also referred as the signature operator \cite[I, p63]{aps}. 

For $s \in {\Bbb C}, \: \: \re(s) > > 0$, put
\[ \eta_A(s) = \frac{2}{\sqrt{\pi}}
\int_0^\infty t^s \: \Tr \left(A e^{- tA^2} \right) dt. \]
Then $\eta_A(s)$ extends to a meromorphic function on ${\Bbb C}$ which
is holomorphic near $s = 0$.
The eta-invariant of $A$ is defined by $\eta(A) = \eta_A(0)$.

\Rk As discussed in the previous remark, the eta invariant of $A$ is twice of the eta invariant
of the signature operator in \cite[I, p63]{aps}.

Let $W$ be an even dimensional compact oriented manifold with boundary $\pt W=M$. If $h$ is a Riemannian metric on $W$ which is of product type near the boundary, then the celebrated Atiyah-Patodi-Singer index formula gives
\be \label{apsit}
 \sign(W) =\int_W {\rm L}(\o{R}{4\pi}) - \hf \eta(A)
\ee
where L denotes the L-polynomial and $R$ the curvature tensor of $h$, and $A$ is the signature operator on $M$ associated with $g=h|_M$. 

Although one usually requires the product metric structure near the boundary, the restriction can be relaxed to that of the boundary being totally geodesic. This is because, in the general case where metric may not be product near the boundary, the extra term coming in is the integral of a transgression of the L-polynomial which depends on the second fundamental form and can be shown to vanish in the totally geodesic case. More precisely, let $h$ be an arbitrary metric on $W$ and $g=h|_M$. If we denote by $x$ the geodesic distance to $M$ in a collar neighborhood, then by the Gauss' lemma, we have
\[ h=dx^2 + g(x) \]
in a collar neighborhood, say $[0, 1]\times M$, and if we denote $y$ to be the local coordinates on $M$, $g(x)=a_{ij}(x, y) dy_i dy_j$. One easily verifies that
\be \label{gisff}
a_{ij}(x, y)= a_{ij}(y) - 2x S_{ij} + O(x^2)
\ee
where $S_{ij}= <\! A(\o{\pt}{\pt y_i}, \ \o{\pt}{\pt y_j}), \o{\pt}{\pt x}\! >$. Now we extend $h$ to a metric $\bar{h}$ on $\overline{W}=[-1, 0] \times M \cup W$ so that, on $[-1, -\hf] \times M$, $\bar{h} = dx^2 + g$ (say, by extending the $a_{ij}(x, y)$).
Applying the index formula of Atiyah-Patodi-Singer to $\overline{W}$ yields
\begin{eqnarray*}
 \sign(W) & = & \sign(\overline{W})= \int_{\overline{W}} {\rm L}(\o{R}{4\pi}) - \hf \eta(A) \\
& = &  \int_{W} {\rm L}(\o{R}{4\pi})  + \int_{[-1, 0] \times M} {\rm L}(\o{R}{4\pi})- \hf \eta(A)
\end{eqnarray*} 
Now, on $[-1, 0] \times M$, $L=dQ$ where $Q$ is the transgression of $L$. More specificly, if we denote by $h_0 =  dx^2 + g$ the product metric on $[-1, 0] \times M$,  $\omega$, $\omega_0$ the connection $1$-forms of $h$, $h_0$ respectively (reverting back to the Cartan formalism) and  $\Omega$, $\Omega_0$ the corresponding curvature $2$-forms, then
\[ Q=r \int_0^1 {\rm L} (\omega -\omega_0, \Omega_t, \dots, \Omega_t) dt, \]
where $r=[\o{n}{4}]$ is the degree of the L-polynomial and $\Omega_t= d\omega_t + \omega_t\wedge \omega_t$ is the curvature form of $\omega_t= t \omega + (1-t) \omega_0$. Thus
\[ \int_{[-1, 0] \times M} {\rm L} = - \int_M Q|_{x=0} .\]
On the other hand, at $x=0$, one verifies that 
\[ \theta = \omega -\omega_0 =\left( \begin{array}{cccc} 0 & \alpha_1 & \cdots & \alpha_n \\
                                      -\alpha_1 & \beta_{11} & \cdots & \beta_{1n}   \\
                                                     \vdots &  \vdots &        &   \vdots \\
                                     -\alpha_n & \beta_{n1} & \cdots & \beta_{nn}
 \end{array} \right) \]
where $\alpha_i = - S(\o{\pt}{\pt y_i}, \o{\pt}{\pt y_j}) dy_j$ and $\beta_{ij} = b_{ij} dx$. Thus, $Q|_M$ depends only on the second fundamental form. Since L is homogenious, when $M$ is totally geodesic, we have 
\[ \int_M Q|_{x=0} \]
as claimed.

Altenatively one may use (\ref{gisff}) to see that the metric is approximately a product up to the second order when $M$ is totally geodesic and argue that the original argument of Atiyah-Patodi-Singer goes through in this case, see the remark after Theorem 4.2 in \cite[I]{aps}.

We recall the following fact from \cite[II, pp 420-421]{aps}. 

\begin{lem} The eta invariant $\eta(A)$ is a conformal invariant. That is, if $\bg = e^{2\phi} g$
is a metric conformal to $g$, and if we denote by $A_g$ and $A_{\bg}$ the corresponding signature operators, we have
\be \label{cioei}
\eta(A_g) = \eta(A_{\bg}).
\ee
\end{lem}

\Pf For completeness and also for comparing with Theorem \ref{co}, we include the following proof from \cite[II]{aps}. Choose a function $a(t)$ such that $a(t) =0$ for small $t$ and $a(t) = 1$ near $t=1$. Consider  the metric $h= e^{2a(t) \phi} (dt^2 + g)$ on $W=M\times [0, \ 1]$, which is of product type near the boundaries. Therefore, by the Atiyah-Patodi-Singer index formula, we have
\[ 0=\sign(W) =\int_W {\rm L}(\o{R}{4\pi}) + \hf (\eta(A_{\bg}- \eta(A_g)). \]
Now, since the L-polynomial is expressible in terms of the Pontryagin forms
which depends only on the Weyl curvature tensor,  the differential form ${\rm L}(\o{R}{4\pi})$ is conformally invariant. But $h$ is conformal to the product metric, giving ${\rm L}(\o{R}{4\pi})=0$. \qed

Thus, since the eta invariant is independent of the choice of the metric in the conformal class, we will denote it by $\eta(M, [g])$. It is a function on the space ${\cal C\! F}(n)$ which is smooth except for simple jump discontinuity on the variation of the metric $g$.

Similarly we have
\begin{theo} \label{co} If $(M, [g])$ is conformally cobordant to zero, then the eta invariant of its signature operator is an even integer. 
\end{theo}

\Pf By definition, if  $(M, [g])$ is conformally cobordant to zero, then it bounds a conformally 
flat manifold $(W, [h])$ such that $[h|_M]=[g]$ and $M$ is a totally umbilic hypersurface of $W$.
According to Lemma \ref{tu}, there is a metric in the conformal class $[h]$, say $h$, such that $M$ is totally geodesic with respect to $h$. Therefore we can apply the Atiyah-Patodi-Singer index formula and obtain
\[ \sign(W) =\int_W {\rm L}(\o{R}{4\pi}) - \hf \eta(M, [g])= - \hf \eta(M, [g]) \]
again by the conformal invariance of the Pontryagin forms. \qed

Thus, as in \cite{lr}, the eta invariant is an obstruction for conformally flat structures bounding in the conformal category. We also note that the eta invariant is reasonably well behaved in the connect sum.

\begin{prop} The function defined on ${\cal C\! F}(n)$ by
\[ \Phi(M, [g])=e^{\pi \eta(M, [g])} \]
is multiplicative under connected sum operation:
\[ \Phi(M_1 \# M_2)= \Phi(M_1) \Phi(M_2). \]
\end{prop}

\Pf First of all, by the conformal flatness of $M_j$, through a conformal deformation of metrics, one can make the metrics on $M_j$ to be the Euclidean metric near any particular point. Let $D_j(3r)\subset M_j$ be round disks of radius $3r$ and $M_{0j}=M_j- ({\rm int} D_j(2r))$. Since the neck region $N_j(3r, r)= D_j(3r) - ({\rm int} D_j(r))$ is conformally equivalent to a cylinder, by a conformal deformation of metrics, one can make the metrics on $M_j$ to be product on $N_j(3r, r)$. Then $ M_1 \# M_2= M_{01} \cup M_{02}$. Now we use the gluing law of eta invariant \cite{df} to deduce
\[ e^{\pi \eta(M_1 \# M_2)} = e^{\pi \eta(M_{01})} e^{\pi \eta(M_{02})} , \]
\[ e^{\pi \eta(M_j)} = e^{\pi \eta(M_{0j})} e^{\pi \eta(D_j(2r))} , \]
\[ e^{\pi \eta(D_1(2r))} e^{\pi \eta(D_2(2r))} = e^{\pi \eta({\rm S}^n)} , \]
where $\eta({\rm S}^n)$ is the eta invariant of the standard sphere which is zero. Combining these equations we obtain the desired result. \qed
\newline

Note that the function $\Phi$ restricts to a group homomorphism from the conformal cobordism group of hyperbolic structures to the circle group
\[ \Phi:\ \Omega_n^{CF}(H) \longrightarrow S^1. \]
 
As pointed in \cite{mn} (Cf. also \cite{lr}) the eta invariant of hyperbolic $3$-manifolds takes values in a dense set of the real line. Thus, the homomorphism $\Phi$ can be highly nontrivial.

X. D: Department of Mathematics, University of  California,
Santa Barbara, California 93106, USA

{\it E-mail: dai@math.ucsb.edu}


\begin{thebibliography}{DSW}
\bibitem[APS]{aps} M.F.Atiyah, V.K.Patodi and I.M.Singer, {\em Spectral asymmetry 
and Riemannian
geometry} I., II., III., Proc. Cambridge Philos. Soc. 77(1975), 78(1975), 79(1976), pp 43-69, pp 405-432,  pp 71-99. 
\bibitem[DF]{df} X. Dai, D. Freed, {\em $\eta$-invariants and determinant lines}
 J. Math. Phys. 35 (1994), 5155-5194.

\bibitem[KP]{kp} R. Kulkarni, U. Pinkall, eds., {\em Conformal Geometry}, Aspects of Mathematics, Vieweg, 1988.
\bibitem[LR]{lr} D. Long, A. Reid, {\em On the geometric boundaries of hyperbolic 4-manifolds}, 
Geometry and Topology, 4(2000), 171-178. 
\bibitem[MN]{mn} R. Meyerhoff, W. Neumann, {\em An aymptotic formula for the eta invariant for hyperbolic 3-manifolds}, Comment. Helv. Math. 67 (1992), pp 28-46.
\end{thebibliography}
\end{document}